\documentclass[reqno]{amsart}
\usepackage{latexsym}
\usepackage{amscd,amsthm,amsfonts,amssymb,amsmath}
\usepackage[all]{xy}
\usepackage[colorlinks=true]{hyperref}
\usepackage{mathrsfs}
\usepackage{stmaryrd}

\newcommand\bs{\backslash}
\newcommand\ra{\rightarrow}
\newcommand\lra{\longrightarrow}
\newcommand\ov{\overline}
\newcommand\ldb{\llbracket}
\newcommand\rdb{\rrbracket}
\newcommand\set[1]{\left\{#1\right\}}
\newcommand\Arr[2]{\left(#1\right)_{\left[#2\right]}}
\newcommand\bfp{{\mathbf{p}}}
\newcommand\bfu{{\mathbf{u}}}
\newcommand\bfx{{\mathbf{x}}}
\newcommand\bfy{{\mathbf{y}}}
\newcommand\bfC{{\mathbf{C}}}
\newcommand\BC{{\mathbb{C}}}
\newcommand\BF{{\mathbb{F}}}
\newcommand\BN{{\mathbb{N}}}
\newcommand\BQ{{\mathbb{Q}}}

\newcommand\BZ{{\mathbb{Z}}}
\newcommand\CA{{\mathcal{A}}}
\newcommand\CB{{\mathcal{B}}}
\newcommand\CL{{\mathcal{L}}}
\newcommand\CN{{\mathcal{N}}}
\newcommand\CR{{\mathcal{R}}}
\newcommand\Gal{{\mathrm{Gal}}}
\newcommand\NP{{\mathrm{NP}}}
\newcommand\Nm{{\mathrm{Nm}}}
\newcommand\ord{{\mathrm{ord}}}
\newcommand\sgn{{\mathrm{sgn}}}
\newcommand\Tr{{\mathrm{Tr}}}

\newtheorem{theorem}{Theorem}[section]
\newtheorem{corollary}[theorem]{Corollary}
\newtheorem{lemma}[theorem]{Lemma}
\newtheorem{proposition}[theorem]{Proposition}
\newtheorem{conjecture}[theorem]{Conjecture}
\theoremstyle{definition}

\newtheorem{example}[theorem]{Example}
\theoremstyle{remark}

\numberwithin{equation}{section}

\begin{document}
\title{On the Newton polygons of twisted $L$-functions of binomials}
\author{Shenxing Zhang}
\date{\today}
\address{School of Mathematics, Hefei University of Technology, Hefei, Anhui 230009, China}
\email{zsxqq@mail.ustc.edu.cn}

\begin{abstract}
Let $\chi$ be an order $c$ multiplicative character of a finite field and $f(x)=x^d+\lambda x^e$ a binomial with $(d,e)=1$.
We study the twisted classical and $T$-adic Newton polygons of $f$.
When $p>(d-e)(2d-1)$, we give a lower bound of Newton polygons and show that they coincide if $p$ does not divide a certain integral constant depending on $p\bmod {cd}$.

We conjecture that this condition holds if $p$ is large enough with respect to $c,d$ by combining all known results and the conjecture given by Zhang-Niu.
As an example, we show that it holds for $e=d-1$. 
\end{abstract}

\maketitle
\tableofcontents

\section{Introduction}
\subsection{Background}
\label{ssec:backg}

Fix a rational prime $p$.
For $q=p^a$ a power of $p$, denote by $\BF_q$ the finite field with $q$ elements, $\BQ_q$ the unramified extension of $\BQ_p$ of degree $a$ and $\BZ_q$ its ring of integers.
Let $f(x)\in\BF_q[x]$ be a polynomial of degree $d$ with Teichm\"uller lifting $\hat f(x)\in\BZ_q[x]$.
Let $\chi:\BF_q^\times\to \BC_p^\times$ be a multiplicative character and $\omega:\BF_q^\times\to \BZ_q^\times$ the Teichm\"uller lifting.
Then we can write $\chi=\omega^{-u}$ for some $0\le u\le q-2$.

For a non-trivial additive character $\psi_m:\BZ_p\to \BC_p^\times$ of order $p^m$, define the twisted $L$-function
	\[L_u(s,f,\psi_m)=\exp\left(\sum_{k=1}^\infty S_{k,u}(f,\psi_m)\frac{s^m}m\right),\]
where $S_{k,u}(f,\psi_m)$ is the twisted exponential sum
	\[S_{k,u}(f,\psi_m)=\sum_{x\in\BF_{q^k}^\times}\psi_m\left(\Tr_{\BQ_{q^k}/\BQ_p}\bigl(\hat f(\hat x)\bigr)\right)\omega^{-u}\left(\Nm_{\BF_{q^k}/\BF_q}(x)\right).\]
If $p\nmid d$, then $L_u(s,f,\psi_m)$ is a polynomial of degree $p^{m-1}d$ by Adolphson-Sperber \cite{AdolphsonSperber1987, AdolphsonSperber1991, AdolphsonSperber1993}, Li \cite{Li1999}, Liu-Wei \cite{LiuWei2007} and Liu \cite{Liu2007}.

We will use the twisted $T$-adic exponential sums developed by Liu-Wan \cite{LiuWan2009} and Liu \cite{Liu2002, Liu2009}.
Define the twisted $T$-adic $L$-function
	\[L_u(s,f,T)=\exp\left(\sum_{k=1}^\infty S_{k,u}(f,T)\frac{s^k}{k}\right)\in 1+s\BZ_q\ldb T\rdb\ldb s\rdb,	\]
where $S_{k,u}(f,T)$ is the twisted $T$-adic exponential sum
	\[S_{k,u}(f,T)=\sum_{x\in\BF_{q^k}^\times}(1+T)^{\Tr_{\BQ_{q^k}/\BQ_p}(\hat f(\hat x))}\omega^{-u}\left(\Nm_{\BF_{q^k}/\BF_q}(x)\right).\]
Then $L_u(s,f,\psi_m)=L_u(s,f,\pi_m)$ where $\pi_m=\psi_m(1)-1$.

Denote by
	\[C_u(s,f,T)=\prod_{j=0}^\infty L_u(q^j s,f,T)\in 1+s\BZ_q\ldb T\rdb\ldb s\rdb\]
the characteristic function, which is $T$-adic entire in $s$.
Then
	\[L_u(s,f,T)=C_u(s,f,T)C_u(qs,f,T)^{-1}.\]
Since the $\pi_m^{a(p-1)}$-adic Newton polygon of $C_u(s,f,\pi_m)$ does not depend on the choice of $\psi_m$, we denote it by $\NP_{u,m}(f)$.
Denote by $\NP_{u,T}(f)$ the $T^{a(p-1)}$-adic Newton polygon of $C_u(s,f,T)$.
As shown in \cite{LiuWan2009} and \cite{Liu2007}, $\NP_{u,m}(f)$ lies over the infinity $u$-twisted Hodge polygon $H_{[0,d],u}^\infty$, which has slopes
\begin{equation}\label{eq:hodge}
	\frac nd+\frac{1}{bd(p-1)}\sum_{k=1}^b u_k,\ n\in\BN.
\end{equation}
If we write $0\le s_0\le\cdots\le s_{p^{m-1}d-1}\le 1$ the $q$-adic slopes of $L_u(s,f,\pi_m)$, then the $q$-adic slopes of $C_u(s,f,\pi_m)$ are
	\[j+s_i,\quad 0\le i\le p^{m-1}d-1, j\in\BN.\]
That's to say, the $\pi_m^{a(p-1)}$-adic Newton polygon of $L_u(s,f,\pi_m)$ is the restriction of $\NP_{u,m}(f)$ on $[0,p^{m-1}d]$, and it determines $\NP_{u,m}(f)$.

The prime $p$ is required large enough in the following results.
When $\chi=\omega^{-u}$ is trivial, in \cite{Zhu2014} and \cite{LiuLiuNiu2009}, they gave a lower bound of the Newton polygons.
They defined a polynomial on the coefficients of $f$, called Hasse polynomial. If the Hasse polynomial is nonzero, then the Newton polygons coincide this lower bound.

Assume that $f(x)=x^d+\lambda x^e$ is a binomial. Since the exponential sums can be transformed to the twisted case when $d$ and $e$ are not coprime, we assume $(d,e)=1$ in this paper.
When $u=0$, we list the known cases here.
\begin{itemize}
\item $p\equiv 1\bmod d$, it's well-known that the Newton polygons coincides the Hodge polygon.
\item $e=1$, see \cite[\S 1, Theorem]{Yang2003}, \cite[Theorem~1.1]{Zhu2014} and \cite[Theorem~1.1]{OuyangYang2016a}.
\item $e=d-1,p\equiv -1\bmod d$, see \cite{OuyangZhang2016}.
\item $e=2,p\equiv 2\bmod d$, see \cite{ZhangNiu2021}. 
\end{itemize}
For arbitrary $u$, Liu-Niu \cite{LiuNiu2011} obtained the Newton polygons when $e=1$. Zhang-Niu \cite{ZhangNiu2021} also give a conjectural description of the Newton polygons when $p\equiv e\bmod d$.

\subsection{Notations}
We list the notations we will use.
\begin{itemize}
\item $i,j,v,w,k,\ell,n$ indices.
\item $f(x)=x^d+\lambda x^e\in \BF_q[x]$ a binomial with $d>e\ge 1,(d,e)=1,\lambda\neq 0$.
\item $\omega^{-u}:\BF_q^\times\to\BC_p^\times$, where $\omega$ is the Teichm\"ller lifting and $0\le u\le q-2$.
\item $H_{[0,d],u}^\infty$, the infinity $u$-twisted Hodge polygon with slopes in \eqref{eq:hodge}.
\item $c=\frac{q-1}{(q-1,u)}$ the order of $\omega^{-u}$, then $u=\frac{(q-1)\mu}c$ for some $(\mu,c)=1$.
\item $P_{u,e,d}$ a polygon with slopes $w(i)$, defined in \eqref{eq:Pchiped}.
\item $b$ the least positive integer such that $p^bu\equiv u\bmod{(q-1)}$ (equivalently, $p^b\equiv 1\bmod c$).
\item $0\le u_i\le p-1$ such that $u=u_0+u_1p+\cdots+u_{a-1}p^{a-1}$, $u_i=u_{b+i}$.
\item $\ov x$ the minimal non-negative residue of $x$ modulo $d$.
\item $\delta_P$ takes value $1$ if $P$ happens; $0$ if $P$ does not happen.
\item $I_n=\set{1,\dots,n}, I_n^*=\set{0,1,\dots,n}$.
\item $S_n$ (resp. $S_n^*$) the set of permutations of $I_n$ (resp. $I_n^*$).
\item $C_{t,n}$ the minimum of $\sum_{i=0}^n \ov{e^{-1}(pi-\tau(i)+t)}$ for $\tau\in S_n^*$ and $S_{t,n}^\circ$ the set of $\tau\in S_n^*$ such that the summation reaches minimal. Set $C_{t,-1}=0$ for convention.
\item $R_{i,\alpha}=\ov{e^{-1}(pi+\alpha)},\ r_{i,\alpha}=\ov{e^{-1}(t-\alpha-i)}$, see Proposition~\ref{pro:polygon}. We will the subscript $\alpha$ if there is no confusion.
\item $\bfC_{t,n,\alpha}$ the maximal size of $\set{i\in I_n^*\mid R_{i,\alpha}+r_{\tau(i),\alpha}\ge d}$ for $\tau\in S_n^*$. We will the subscript $\alpha$ if there is no confusion.

\item $y_{t,i}^\tau=\ov{e^{-1}(pi-\tau(i)+t)},\ x_{t,i}^\tau=d^{-1}(pi-\tau(i)+t-ey_{t,i}^\tau)$ the unique solution of $dx+ey=pi-\tau(i)+t$ with $0\le y\le d-1$.
\item $h_{n,k}, h_{u,e,d}$ the Hasse numbers defined in \eqref{eq:hasseh}.
\item $\bfp$ the minimal non-negative residue of $p$ modulo $cd$.
\item $H_{\mu,c,\bfp,e,d}\in\BZ$ a constant defined in \eqref{eq:hasse}.

\item $E(X)$ the $p$-adic Artin-Hasse series, see \eqref{eq:artin-hasse}.
\item $\pi$ a $T$-adic uniformizer of $\BQ_p\ldb T\rdb$ given by $E(\pi)=1+T$, with a fixed $d(q-1)$-th root $\pi^{\frac1{d(q-1)}}$.
\item $E_f(X)$, see \eqref{eq:E_f}.
\item $M_u=\frac{u}{q-1}+\BN$.
\item $\CL_u$ a Banach space, see \eqref{eq:CL_u}.
\item $\CB_u$ a subspace of $\CL_u$, see \eqref{eq:CB_u}.
\item $\CB=\CB_u\oplus\CB_{pu}\oplus\cdots\oplus\CB_{p^{b-1}u}$.
\item $\psi:\CL_u\ra \CL_{p^{-1}u}$ defined as $\psi\left(\sum_{v\in M_u} b_v X^v\right)=\sum_{v\in M_{p^{-1}u}} b_{pv}X^v$.
\item $\sigma\in\Gal(\BQ_q/\BQ_p)$ the Frobenius, which acts on $\CL_u$ via the coefficients.
\item $\Psi=\sigma^{-1}\circ\psi\circ E_f:\CB_u\to\CB_{p^{-1}u}$ the Dwork's $T$-adic semi-linear operator.
\item $c_n$ the coefficients of $\det(1-\Psi s\mid \CB)$, see \eqref{eq:c_n}.

\item $s_k\equiv p^k u\bmod{q-1}$ with $0\le s_k\le q-2$.
\item $\Gamma=\left(\gamma_{(v,\frac{s_k}{q-1}+i),(w,\frac{s_\ell}{q-1}+j)}\right)$ the matrix coefficient of $\Psi$ on $\CB$, see  \eqref{eq:Gamma}.
\item $\Gamma^{(k)}$ the sub-matrix of $\Gamma$ defined in \eqref{eq:Gamma}.
\item $A^{(k)}=A\cap\Gamma^{(k)}$ the sub-matrix of a principal minor $A$ of $\Gamma$.
\item $\CA_n$ the set of all principal minor $A$ of order $bn$, such that every $A^{(k)}$ has order $n$.
\item $\phi(n)\in\BN\cup\set{+\infty}$ the minimal $x+y$ where $dx+ey=n,x,y\in\BN$.
\item $\gamma_{(\frac{s_k}{q-1}+i,\frac{s_\ell}{q-1}+j)}$, see \eqref{eq:gamma2}.
\item $\Arr{x}{n}:=x(x-1)\cdots(x-n+1), \Arr{x}{0}:=1$ the falling factorial.
\end{itemize}

\subsection{Main results}
In this paper, we give an explicit lower bound of Newton polygons of twisted $L$-functions of binomial $f(x)=x^d+\lambda x^e$.
We reduce the Hasse polynomial to a certain integer \eqref{eq:hasse}. Then $p>(d-e)(2d-1)$ does not divide this constant, if and only if this lower bound coincides the Newton polygons.
Finally, we show that this condition holds for $e=d-1$.

Denote by $P_{u,e,d}$ the polygon such that
\begin{equation}\label{eq:Pchiped}
	P_{u,e,d}(n)=\frac{n(n-1)}{2d}+\frac{1}{bd(p-1)}\sum_{k=1}^b \bigl(nu_k+(d-e)C_{u_k,n-1}\bigr),\ n\in\BN.
\end{equation}
Denote by $w(n)=P_{u,e,d}(n+1)-P_{u,e,d}(n)$. Then
	\[w(n)=\frac nd+\frac{1}{bd(p-1)}\sum_{k=1}^b\bigl(u_k+(d-e)(C_{u_k,n}-C_{u_k,n-1})\bigr).\]
This polygon lies above the Hodge polygon $H_{[0,d],u}^\infty$ with same points at $d\BZ$, and $w(n+d)=1+w(n)$. Moreover, we have $w(n)\le w(n+1)$ if $p>(d-e)(2d-1)$. See Proposition~\ref{pro:polygon}.

\begin{theorem}\label{thm:lower_bound}
Assume that $p>(d-e)(2d-1)$. Then $\NP_{u,T}(f)$ lies above $P_{u,e,d}$. As a corollary, $\NP_{u,m}(f)$ lies above $P_{u,e,d}$.
\end{theorem}

Define
\begin{equation}\label{eq:hasseh}
	h_{n,k}:=\sum_{\tau\in S_{u_k,n}^\circ}\sgn(\tau)\prod_{i=0}^n\frac{1}{x_{u_k,i}^\tau!y_{u_k,i}^\tau!},\quad	h_{u,e,d}:=\prod_{n=0}^{d-2}\prod_{k=1}^b h_{n,k}.
\end{equation}

\begin{theorem}\label{thm:coincide}
Assume that $p>(d-e)(2d-1)$. Then
\begin{equation}\label{eq:equal}
	\NP_{u,m}(f)=\NP_{u,T}(f)=P_{u,e,d}
\end{equation}
holds if and only if $h_{u,e,d}\in\BZ_p^\times$, if and only if $p\nmid H_{\mu,c,\bfp,e,d}$.
\end{theorem}

Here $H_{\mu,c,\bfp,e,d}\in\BZ$ is a constant defined in \eqref{eq:hasse} and $\bfp$ is the minimal positive residue of $p$ modulo $cd$. Thus we have the following corollary.
\begin{corollary}\label{cor:cong}
Assume that \eqref{eq:equal} holds for
	\[a, m, p, f(x)=x^d+\lambda x^e\in\BF_{p^a}[x], u=\frac{(p^a-1)\mu}c,\]
where $b\mid a, \lambda\neq 0$ and $p>(d-e)(2d-1)$.
Then
\begin{enumerate}
\item $H_{\mu,c,\bfp,e,d}\neq 0$.
\item For any
	\[a', m', p',f'(x)=x^d+\lambda' x^e\in\BF_{p^{\prime a'}}[x], u'=\frac{({p'}^{a'}-1)\mu}c,\]
where $b\mid a, \lambda\neq0$ and $p'>(d-e)(2d-1)$, we have \eqref{eq:equal} if $p'\equiv p\bmod {cd}$ and $p'>H_{\mu,c,\bfp,e,d}$.
\item As $p'\equiv p\bmod cd$ tends to infinity, the polygons $\NP_{u,m}(f)$ and $\NP_{u,T}(f)$ tend to $H_{[0,d],u}^\infty$, which only depends on $\mu,c,\bfp,d$.
\end{enumerate}
\end{corollary}

The following result extends \cite{OuyangZhang2016}, as they considered the untwisted case with an additional condition $p\equiv -1\bmod d$.

\begin{theorem}\label{thm:examples}
Assume that $e=d-1$.
We have $\NP_{u,m}(f)=\NP_{u,T}(f)=P_{u,e,d}$ if $p>c(d^2-d+1)$.
\end{theorem}

We give the following conjecture, which generalizes the conjecture in \cite{ZhangNiu2021}.
Note that $h_{u,e,d}$ may be zero since $S_{u_k,n}^\circ$ may be empty, so we require that $p$ is large with respect to $c$, as in Corollary~\ref{cor:cong} and Theorem~\ref{thm:examples}.
\begin{conjecture}\label{con:lower_bound}
If $p$ is large enough with respect to $c,d$, then $\NP_{u,m}(f)=\NP_{u,T}(f)=P_{u,e,d}$.
\end{conjecture}

\section{The lower bound}

\subsection{The property of the lower bound polygon}

For any integer $t$, we denote
	\[C_{t,n}=\min_{\tau\in S_n^*} \sum_{i=0}^n \ov{e^{-1}(pi-\tau(i)+t)}.\]
We set $C_{t,-1}=0$ for convention.
For any integer $\alpha$, we denote
	\[R_{i,\alpha}=\ov{e^{-1}(pi+\alpha)},\ r_{i,\alpha}=\ov{e^{-1}(t-\alpha-i)}\]
and
	\[\bfC_{t,n,\alpha}=\max\#\set{i\in I_n^*\mid R_{i,\alpha}+r_{\tau(i),\alpha}\ge d}.\]

\begin{proposition}\label{pro:polygon}
(1) For any $\alpha$, we have
	\[C_{t,n}=\sum_{i=0}^n(R_{i,\alpha}+r_{i,\alpha})-d \bfC_{t,n,\alpha}.\]

(2) For any $\alpha$, we have
	\[\bfC_{t,n+d,\alpha}=d-1+\bfC_{t,n,\alpha},\quad C_{t,n+d}=C_{t,n}.\]
Thus $w(n+d)=1+w(n)$ and $P_{u,e,d}(dn)=H_{[0,d],u}^\infty(dn)$.

(3) If $p>(d-e)(2d-1)$, we have $w(n)\le w(n+1)$.
\end{proposition}

\begin{proof}
We omit the subscript $\alpha$ in this proof for convention.

(1) It follows from
	\[\ov{e^{-1}(pi-\tau(i)+t)}=R_i+r_{\tau(i)}-d\delta_{R_i+r_{\tau(i)}\ge d}.\]
	
(2) We have 
	\[\bfC_{t,n}=\max_{\tau\in S_n^*}\#\set{i\in I_n^*\mid R_i\ge d-r_{\tau(i)}}.\]
Note that 
	\[\set{R_i\mid i\in I_{n+d}^*}=\set{R_i\mid i\in I_n^*}\cup\set{0,1,\dots,d-1}, \]
	\[\set{d-r_i\mid i\in I_{n+d}^*}=\set{d-r_i\mid i\in I_n^*}\cup\set{d,1,\dots,d-1}. \]
We may drop the $0$ and $d$ since they do not affect the size.
Apple Lemma~\ref{lem:remove_one_element} $(d-1)$ times and we get $\bfC_{t,n+d}=d-1+\bfC_{t,n}$.

Since
	\[\sum_{i=n+1}^{n+d}(R_i+r_i)=2\sum_{j=0}^{d-1}j=d(d-1),\]
we have $C_{t,n+d}=C_{t,n}$. Thus $w(n+d)=1+w(n)$.

Note that $C_{t,n+d}=C_{t,n}$ also holds for $n=-1$. Hence $C_{t,dn-1}=0$ and $P_{u,e,d}(dn)=H_{[0,d],u}^\infty(dn)$.

(3) Denote by $\delta=\delta_{R_n+r_n\ge d}$. For any $\tau\in S_n^*$, write $i=\tau(n)$, $j=\tau^{-1}(n)$ and $\tau_1=(ni)\tau$. Then $\tau_1(n)=n$, $\tau_1(j)=i$ and
	\[\begin{split}
	&\delta+\#\set{i\in I_{n-1}^*\mid R_i+r_{\tau_1(i)}\ge d}-\#\set{i\in I_n^*\mid R_i+r_{\tau(i)}\ge d}\\
	=&\delta+\delta_{R_j+r_i\ge d}-\delta_{R_j+r_n\ge d}-\delta_{R_n+r_i\ge d}.
	\end{split}\]
If this is $-2$, then $2d>R_n+r_n+R_j+r_i\ge 2d$, that's impossible. Thus $\delta+\bfC_{t,n-1}-\bfC_{t,n}\ge -1$.

Any $\sigma\in S_{n-1}^*$ can be viewed as an element $\sigma_1\in S_n^*$ fixing $n$. Thus
	\[\delta+\#\set{i\in I_{n-1}^*\mid R_i+r_{\sigma(i)}\ge d}=\#\set{i\in I_n^*\mid R_i+r_{\sigma_1(i)}\ge d}.\]
and then $\delta+\bfC_{t,n-1}\le \bfC_{t,n}$.

Now
	\[\begin{split}
	 &C_{t,n}-C_{t,n-1}\\
	=&R_n+r_n-d(\bfC_{t,n}-\bfC_{t,n-1})\\
	=&\ov{e^{-1}(pn-n+t)}+d(\delta+\bfC_{t,n-1}-\bfC_{t,n})
	\end{split}\]
lies in $[-d,d-1]$.
Therefore,
	\[\begin{split}
	 &w(n)-w(n-1)\\
	=&\frac{1}{d}+\frac{d-e}{bd(p-1)}\sum_{k=1}^b(C_{u_k,n}-2C_{u_k,n-1}+C_{u_k,n-2})\\
	\ge&\frac{1}{d}+\frac{(d-e)(1-2d)}{d(p-1)}\ge 0
	\end{split}\]
since $p>(d-e)(2d-1)$.
\end{proof}

\begin{lemma}\label{lem:remove_one_element}
Let $A=\set{a_0,\dots,a_m}$ and $B=\set{b_0,\dots,b_m}$ be two multi-sets of integers.
Assume that $a_0\ge b_0$ and for any $i>0$, $b_i>a_0$ or $b_i\le b_0$.
Then
	\[\max_{\tau\in S_m^*}\#\set{i\in I_m^*\mid a_i\ge b_{\tau(i)}}=1+\max_{\sigma\in S_m}\#\set{i\in I_m\mid a_i\ge b_{\tau(i)}}.\]
\end{lemma}
\begin{proof}
Every permutation in $S_n$ can be viewed as an permutation in $S_n^*$ fixing $0$, thus ``$\ge$'' holds trivially. Write $i=\tau(0)$, $j=\tau^{-1}(0)$ and $\tau_1=(0i)\tau$. Then $\tau_1(0)=0$ and $\tau_1(j)=i$. Thus
	\[\begin{split}
	&\#\set{i\in I_m^*\mid a_i\ge b_{\tau_1(i)}}-\#\set{i\in I_m^*\mid a_i\ge b_{\tau(i)}}\\
	=&1+\delta_{a_j\ge b_i}-\delta_{a_j\ge b_0}-\delta_{a_0\ge b_i}.
	\end{split}\]
If this is negative, then $a_0\ge b_i>a_j\ge b_0$, which is impossible.
Thus ``$\le$'' holds. 
\end{proof}

\subsection{The twisted \texorpdfstring{$T$}{T}-adic Dwork's trace formula}

This part is almost the same with \cite[\SS 2,3]{LiuNiu2011}.
Denote by
\begin{equation}\label{eq:artin-hasse}
	E(X)=\exp\left(\sum_{i=0}^\infty p^{-i}X^{p^i}\right)=\sum_{n=0}^\infty \lambda_n X^n\in\BZ_p\ldb X\rdb
\end{equation}
the $p$-adic Artin-Hasse series.
Then $\lambda_n=1/n!$ if $n<p$.
Denote by
\begin{equation}\label{eq:E_f}
	E_f(X)=E(\pi X^d)E(\pi\hat\lambda X^e)=\sum_{n=0}^\infty \gamma_n X^n.
\end{equation}
Then
	\[\gamma_k=\sum\pi^{x+y}\lambda_x\lambda_y\hat\lambda^y,\]
where $(x,y)$ runs through non-negative solutions of $dx+ey=k$.

Denote by $M_u=\frac{u}{q-1}+\BN$.
Define
\begin{equation}\label{eq:CL_u}
	\CL_u=\set{\left.\sum_{v\in M_u} b_v \pi^{\frac{v}{d}}X^v\;\right|\; b_v\in\BZ_q\ldb \pi^{\frac{1}{d(q-1)}}\rdb}
\end{equation}
and
\begin{equation}\label{eq:CB_u}
	\CB_u=\set{\left.\sum_{v\in M_u} b_v \pi^{\frac{v}{d}}X^v\in\CL_u\;\right|\;\ord_\pi b_v\to+\infty\text{ as }v\to+\infty}.
\end{equation}
Define a map 
\begin{equation}\label{eq:psi}
	\begin{split}
		\psi:\CL_u&\lra \CL_{p^{-1}u}\\
		\sum_{v\in M_u} b_v X^v&\longmapsto \sum_{v\in M_{p^{-1}u}} b_{pv}X^v.
	\end{split}
\end{equation}
The power series $E_f$ defines a map on $\CB_u$ via multiplication.
Let $\sigma\in \Gal(\BQ_q/\BQ_p)$ be the Frobenius, which acts on $\CL_u$ via the coefficients.
Then the Dwork's $T$-adic semi-linear operator $\Psi=\sigma^{-1}\circ\psi\circ E_f$ sends $\CB_u$ to $\CB_{p^{-1}u}$.
Hence $\Psi$ acts on 
	\[\CB:=\bigoplus_{i=0}^{b-1}\CB_{p^iu}.\]
We have a linear map
	\[\Psi^a=\psi^a\circ\prod_{i=0}^{a-1}E_f^{\sigma^i}(X^{p^i})\]
on $\CB$ over $\BZ_q\ldb\pi^{\frac{1}{d(q-1)}}\rdb$.
Since $\Psi$ is completely continuous in the sense of \cite{Serre1962}, the following determinants are well-defined.

\begin{theorem}\label{thm:characterize Newton polygon}
We have
	\[C_u(s,f,T)=\det\left(1-\Psi^a s \;\left|\; \CB_u/\BZ_q\ldb\pi^{\frac{1}{d(q-1)}}\rdb\right.\right).\]
Thus the $T$-adic Newton polygon of $C_u(s,f,T)$ is the lower convex closure of 
	\[\left(n,\frac1b\ord_T(c_{abn})\right),\ n\in\BN,\]
where
\begin{equation}\label{eq:c_n}
	\det\left(1-\Psi s\;\left|\;\CB/\BZ_p\ldb\pi^{\frac{1}{d(q-1)}}\rdb\right.\right)=\sum_{i=0}^\infty (-1)^n c_n s^n.
\end{equation}
\end{theorem}

\begin{proof}
See \cite[Theorem~4.8]{LiuWan2009}, \cite{Liu2007}, \cite[Theorems~2.1, 2.2]{LiuLiuNiu2009} and \cite[Theorems~2.1, 5.3]{LiuNiu2011}.
\end{proof}

Write $s_k\equiv p^k u\bmod{q-1}$ with $0\le s_k\le q-2$. Then $s_{b-k}=s_{-k}=u_k+u_{k+1}p+\cdots+u_{k+a-1}p^{a-1}$.
Let $\xi_1,\dots,\xi_a$ be a normal basis of $\BQ_q$ over $\BQ_p$.
The space $\CB$ has a basis 
	\[\set{\xi_v(\pi^{\frac1d}X)^{\frac{s_k}{q-1}+i}}_{(i,v,k)\in\BN\times I_a\times I_b}\]
over $\BZ_p\ldb\pi^{\frac{1}{d(q-1)}}\rdb$.
Let $\Gamma=\left(\gamma_{(v,\frac{s_k}{q-1}+i),(w,\frac{s_\ell}{q-1}+j)}\right)_{\BN\times I_a\times I_b}$ be the matrix of $\Psi$ on $\CB$ with respect to this basis.
Then
\begin{equation}\label{eq:Gamma}
	\Gamma=\begin{pmatrix}
		0&\Gamma^{(1)}&0&\cdots&0\\
		0&0&\Gamma^{(2)}&\cdots&0\\
		\vdots&\vdots&\vdots&\ddots&\vdots\\
		0&0&0&\cdots&\Gamma^{(b-1)}\\
		\Gamma^{(b)}&0&0&\cdots&0
	\end{pmatrix},
\end{equation}
where
	\[\Gamma^{(k)}=\left(\gamma_{(v,\frac{s_{k-1}}{q-1}+i),(w,\frac{s_k}{q-1}+j)}\right)_{\BN\times I_a}.\]
Hence we have
	\[\det\left(1-\Psi s\;\left|\; \CB/\BZ_p\ldb\pi^{\frac{1}{d(q-1)}}\rdb\right.\right)=\det(1-\Gamma s)=\sum_{n=0}^\infty (-1)^{bn} c_{bn}s^{bn}\]
with $c_n=\sum \det(A)$, where $A$ runs through all principal minors of order $n$, see \cite{LiZhu2005}.
Denote by $A^{(k)}=A\cap \Gamma^{(k)}$ as a minor of $\Gamma^{(k)}$.
If $A$ has order $bn$, but the order of some $A^{(k)}$ is not $n$, then $\det(A)=0$.
Denote by $\CA_n$ the set of all principal minors of order $bn$, such that every $A^{(k)}$ has order $n$.
Then
\begin{equation}\label{eq:expression in terms of minors}
	c_{bn}=\sum_{A\in\CA_n}\det(A)=(-1)^{n(b-1)}\sum_{A\in\CA_n}\prod_{k=1}^b \det(A^{(k)}).
\end{equation}

\begin{theorem}\label{thm:lower_coe}
If $p>(d-e)(2d-1)$, then
	\[\ord_\pi(\det(A))\ge ab(p-1)P_{u,e,d}(n+1)\]
for any $A\in\CA_{a(n+1)}$.
\end{theorem}

\begin{proof}[Proof of Theorem~\ref{thm:lower_bound}]
By Theorem~\ref{thm:lower_coe} and \eqref{eq:expression in terms of minors}, we have
	\[\ord_\pi(c_{abn})\ge ab(p-1)P_{u,e,d}(n).\]
Thus $\NP_{u,T}(f)$ lies above $P_{u,e,d}$ by Theorem~\ref{thm:characterize Newton polygon}.
Note that $\NP_{u,m}(f)\ge \NP_{u,T}(f)$ by definition. Therefore, $\NP_{u,m}(f)$ also lies above $P_{u,e,d}$.
\end{proof}

\subsection{Estimation on \texorpdfstring{$c_n$}{coefficients}}

Denote by
	\[\phi(n)=\min\set{x+y\mid dx+ey=n,x,y\in\BN}\in\BN\cup\set{+\infty}.\]
Here the minimal element in $\emptyset$ is regarded as $+\infty$.
For $i,j\in\BN,k\in I_b$, define
\begin{equation}\label{eq:gamma2}
	\gamma_{(\frac{s_{k-1}}{q-1}+i,\frac{s_k}{q-1}+j)}=
		\pi^{\frac{s_k-s_{k-1}}{d(q-1)}+\frac{j-i}{d}}\gamma_{pi-j+u_{-k}}.
\end{equation}
Then
	\[\xi_w^{\sigma^{-1}}\gamma_{(\frac{s_{k-1}}{q-1}+i,\frac{s_k}{q-1}+j)}^{\sigma^{-1}}=\sum_{u\in I_a}\gamma_{(v,\frac{s_{k-1}}{q-1}+i),(w,\frac{s_k}{q-1}+j)}\xi_v\]
and
\begin{equation}\label{eq:ord terms}
	\begin{split}
	&\ord_\pi\left(\gamma_{(v,\frac{s_{k-1}}{q-1}+i),(w,\frac{s_k}{q-1}+j)}\right) \ge \ord_\pi\left(\gamma_{(\frac{s_{k-1}}{q-1}+i,\frac{s_k}{q-1}+j)}\right)\\
	=&\frac{s_k-s_{k-1}}{d(q-1)}+\frac{j-i}{d}+\phi(pi-j+u_{-k}).
	\end{split}
\end{equation}

\begin{lemma}\label{lem:lower_bound}
For any $\tau\in S_n^*$ and integer $t$,
	\[\sum_{i=0}^n \phi(pi-\tau(i)+t)\ge d^{-1}\left(\frac{(p-1)n(n+1)}2+(n+1)t+(d-e)C_{t,n}\right).\]
\end{lemma}

\begin{proof}
We may assume that $pi-\tau(i)+t\in d\BN+e\BN$ for each $i$.
One can easily show that
	\[\phi(k)=d^{-1}\left(k+(d-e)\ov{e^{-1}k}\right)\]
and the minimum arrives at
	\[(x,y)=\left(d^{-1}(k-e\ov{e^{-1}k}),\ov{e^{-1}k}\right).\]
Thus
\begin{equation}\label{eq:phi}
	\phi(pi-j+t)=d^{-1}\left(pi-j+t+(d-e) \ov{e^{-1}(pi-j+t)}\right).
\end{equation}
The result then follows easily.
\end{proof}

\begin{lemma}\label{lem:d copies}
Assume $a_i=a_{i+m}$ and $b_i=b_{i+m}$ for any $i$. Then
	\[\max_{\tau\in S_{md}}\#\set{i\in I_{md}\mid a_i\ge b_{\tau(i)}}=d\max_{\sigma\in S_m}\#\set{i\in I_m\mid a_i\ge b_{\sigma(i)}}.\]
\end{lemma}
\begin{proof}
We may assume that $a_k\ge b_k$ and for any $i\neq k$, $b_i>a_k$ or $b_i\le b_k$. Otherwise both sides should be zero. 
We may assume that $k=m$ for simplicity.
Apply Lemma~\ref{lem:remove_one_element} to $(a_{mi},b_{mi})$, we get
	\[\max_{\tau\in S_{md}}\#\set{i\in I_{md}\mid a_i\ge b_{\tau(i)}}
	=d+\max_{\sigma}\#\set{i\in I_{md}-m\BZ\mid a_i\ge b_{\tau(i)}},\]
where $\sigma$ runs through permutations on $I_{md}-m\BZ$.
Since
	\[\max_{\tau\in S_{m}}\#\set{i\in I_{m}\mid a_i\ge b_{\tau(i)}}
	=1+\max_{\sigma}\#\set{i\in I_{m}-\set{m}\mid a_i\ge b_{\tau(i)}},\]
where $\sigma$ runs through permutations on $I_{m}-\set m$, the result then follows by induction on $m$.
\end{proof}

\begin{lemma}\label{lem:lower_bound N}
For any $i\in \BN\times I_a$, we write $i=(i',i'')$.
Then for any permutation $\tau$ on $I_n^*\times I_a$,
	\[\sum_{i\in I_n^*\times I_a} \phi(pi'-\tau(i)'+t)\ge\frac{a}{d}\left(\frac{(p-1)n(n+1)}{2}+(n+1)t+(d-e)C_{t,n}\right).\]
\end{lemma}
\begin{proof}
By Eq.~\eqref{eq:phi}, we only need to show that
 	\[\min_{\tau}\sum_{i\in I_n^*\times I_a}\ov{e^{-1}(pi-\tau(i)+t)}=a C_{t,n}.\]
By Proposition~\ref{pro:polygon}, it can be reduced to
	\[\max_\tau\#\set{i\in I_n^*\times I_a\mid R_{i',\alpha}+r_{\tau(i)',\alpha}\ge d}=a \bfC_{t,n,\alpha}.\]
This follows from Lemma~\ref{lem:d copies}.
\end{proof}

\begin{proof}[Proof of Theorem~\ref{thm:lower_coe}]
This proof is similar to \cite[Theorem~3.2]{ZhangNiu2021}.
Denote by $\CR$ the set of indices of $A$ and
	\[\CR^{(k)}\times\set k=\CR\cap(\BN\times I_a\times\set k),\quad \CR^{(0)}=\CR^{(b)}.\]
Then $\#\CR^{(k)}=a(n+1)$,
	\[A^{(k)}=\left(\gamma_{(v,\frac{s_{k-1}}{q-1}+i),(w,\frac{s_k}{q-1}+j)}\right)_{(i,v)\in \CR^{(k-1)},(j,w)\in\CR^{(k)}}\]
and
	\[\det(A)=\prod_{k=1}^b\det(A^{(k)})=\sum_{\tau}\sgn(\tau)\prod_{i\in\CR}\gamma_{i,\tau(i)},\]
where $\tau$ runs through permutations of $\CR$ such that $\tau(\CR^{(k-1)})=\CR^{(k)}$.
Here,
	\[\ord_\pi\left(\prod_{i\in\CR}\gamma_{i,\tau(i)}\right)\ge S_\CR^\tau\]
by \eqref{eq:ord terms}, where
\[\begin{split}
	S_\CR^\tau&=\sum_{k=1}^b\sum_{i\in \CR^{(k-1)}} \left(\frac{\tau(i)'-i'}{d}+\phi\bigl(pi'-\tau(i)'+u_{-k}\bigr)\right)\\
	&\ge d^{-1}\sum_{k=1}^b\sum_{i\in\CR^{(k-1)}}\left((p-1)i'+(d-e)\ov{e^{-1}(pi'-\tau(i)'+u_{-k})}\right)
\end{split}\]
by Eq.~\eqref{eq:phi}.
By Lemma~\ref{lem:lower_bound N},
	\[S_\CN^\sigma\ge ab(p-1)P_{u,e,d}(n+1),\]
where $\CN=I_n^*\times I_a\times I_b$.
By \eqref{eq:expression in terms of minors}, we only need to show that for any permutation $\tau$ of $\CR\neq\CN$ such that $\tau(\CR^{(k-1)})=\CR^{(k)}$, there is a permutation $\sigma$ of $\CN$ such that $\sigma(\CN^{(k-1)})=\CN^{(k)}$ and $S_\CR^\tau\ge S_\CN^\sigma$. 

Assume $\#(\CR\bs \CN)=m$. Write $T=(\CN\bs \CR)\cup\tau^{-1}(\CR\bs\CN)$, then $\#T=2m$ and $\CN\bs T=\CN\cap\tau^{-1}(\CN\cap\CR)$. Thus $\tau(\CN\bs T)\subset \CN$.
Note that for $i\in\CR\bs\CN,j\in\CN\bs\CR$, $i'\ge n+1\ge j'+1$.
We can choose a permutation $\sigma$ of $\CN$ such that $\sigma(\CN^{(k-1)})=\CN^{(k)}$ and $\sigma=\tau$ on $\CN\bs T$. Then
	\[\begin{split}
&d(S_\CR^\tau-S_\CN^\sigma)\\
\ge&\left(\sum_{i\in \CR\bs\CN}-\sum_{i\in\CN\bs\CR}\right)(p-1)i'-\sum_{k=1}^b\sum_{i\in T\cap\CN^{(k)}}(d-e)\ov{e^{-1}(pi'-\tau(i)'+u_{-k})}\\
\ge&m(p-1)-2m(d-e)(d-1)>0.
\end{split}\]
The result then follows.
\end{proof}

\section{The Newton polygons}

\begin{lemma}\label{lem:rigidity}
The Newton polygon $\NP_m(f)$ lies over $\NP_T(f)$. Moreover, if the equality holds for one $m$, then it holds for all $m$.
\end{lemma}
\begin{proof}
See \cite[Theorem~2.3]{LiuWan2009} and \cite[Theorem~5.5]{LiuNiu2011}.
\end{proof}

\begin{proof}[Proof of Theorem~\ref{thm:coincide}]
(1) Since $w(d+i)=1+w(i)$, both of $\NP_{u,m}(f)$ and $P_{u,e,d}$ across points $\bigl(di,H_{[0,d],u}^\infty(di)\bigr)$, we only need to show that $\NP_{u,m}(f)=P_{u,e,d}$ on $[1,d-1]$. By Lemma~\ref{lem:rigidity}, we may assume that $m=1$.

Assume $0\le n\le d-2$. Recall that $S_{t,n}^\circ$ is the set of $\tau\in S_n^*$ such that
	\[\#\set{i\in I_n^*\mid R_{i,\alpha}+r_{\tau(i),\alpha}\ge d}=\bfC_{t,n,\alpha}\]
and every $pi-\tau(i)+t\in d\BN+e\BN$.
It's equivalently to say, the equality in Lemma~\ref{lem:lower_bound} holds.
Recall that
	\[y_{t,i}^\tau=\ov{e^{-1}(pi-\tau(i)+t)},\quad x_{t,i}^\tau=\phi(pi-\tau(i)+t)-y_{t,i}^\tau.\]
Denote by $m$ the right hand side in Lemma~\ref{lem:lower_bound}.
Then we have
	\[\begin{split}
	&\det(\gamma_{pi-j+t})_{i,j\in I_n^*}
	 \equiv\pi^m\sum_{\tau\in S_{t,n}^\circ}\sgn(\tau)\prod_{i=0}^n\lambda_{x_{t,i}^\tau}\lambda_{y_{t,i}^\tau}\hat\lambda^{y_{t,i}^\tau}
	\\
	\equiv&\pi^m\hat\lambda^{v_{t,n}}\sum_{\tau\in S_{t,n}^\circ}\sgn(\tau)\prod_{i=0}^n\frac{1}{x_{t,i}^\tau!y_{t,i}^\tau!}\mod \pi^{m+1},
	\end{split}\]
where
	\[v_{t,n}:=\sum_{i=0}^n y_{t,i}^\tau=\sum_{i=1}^n(R_{i,\alpha}+r_{i,\alpha})-d\bfC_{t,n,\alpha}\]
is independent on $\tau\in S_n^\circ$.

Recall that $S_\CR^\tau>S_\CN^\sigma$ in the proof of Theorem~\ref{thm:lower_coe}.
Then modulo $\pi^{ab(p-1)P_{u,e,d}(n+1)+1}$, we have
	\[\begin{split}
	c_{ab(n+1)}&=\sum_{A\in\CA_{a(n+1)}}\det(A)\equiv \det\bigl((\gamma_{i,j})_{i,j\in \CN}\bigr)\\
	&=\pm\Nm\left(\prod_{k=1}^b\det\left(\gamma_{(\frac{s_{k-1}}{q-1}+i,\frac{s_k}{q-1}+j)}\right)_{i,j\in I_n^*}\right)\\
	&=\pm \Nm\left(\prod_{k=1}^b\det(\gamma_{pi-j+u_k})_{i,j\in I_n^*}\right)\\
	&\equiv\pm\pi^{ab(p-1)P_{u,e,d}(n+1)}\Nm\left(\prod_{k=1}^b \hat\lambda^{v_{u_k,n}} h_{n,k}\right)
	\end{split}\]
by \eqref{eq:expression in terms of minors}, \eqref{eq:gamma2},  \cite[Lemma~4.4]{LiuLiuNiu2009} and \cite[Lemma~3.5]{LiuNiu2011}.
Hence we get the first assertion by replacing $\pi$ by $\pi_1$.

(2) Denote by $t_k$ the minimal non-negative residue of $p^{-k} \mu$ modulo $c$.
Then $u_k=\frac{t_{k+1} p-t_k}{c}$.
Write $\bfp$ the minimal positive residue of $p$ modulo $cd$ and $p=c d\ell+\bfp$.
Denote by
	\[\bfu_k=\frac{t_{k+1} \bfp-t_k}c,\ 
	\bfy_{\bfu_k,i}^\tau=\ov{-e^{-1}(\bfp i-\tau(i)+\bfu_k)},\ 
	\bfx_{\bfu_k,i}^\tau=\frac{\bfp i-\tau(i)+\bfu_k-e\bfy_{\bfu_k,i}^\tau}{d}.\]
Then
	\[u_k=t_{k+1} d\ell+\bfu_k,\ 
	y_{u_k,i}^\tau=\bfy_{\bfu_k,i}^\tau,\ 
	x_{u_k,i}^\tau=(c i+t_{k+1})\ell+\bfx_{\bfu_k,i}^\tau.\]
It's easy to see that $\bfx_{\bfu_k,i}^\tau<\bfp$ and $x_{u_k,i}^\tau<p$.
Since
	\[\bfx_{\bfu_k,i}^\tau\ge \frac{-n-e(d-1)}{d}>-e-1,\]
we have $\bfx_{\bfu_k,i}^\tau\ge -e$.
Note that $y_{t,i}^\tau$ does not depend on $\ell$.
Denote by 
\begin{equation}\label{eq:hasse}
	\begin{split}
	H_{\mu,c,\bfp,e,d}=&\prod_{k=1}^b\prod_{n=0}^{d-2}\sum_{\tau\in S_n^\circ} \sgn(\tau)\prod_{i=1}^n \Arr{d-1}{d-1-\bfy_{\bfu_k,i}^\tau}
	\times (cd)^{\bfp-1-\bfx_{\bfu_k,i}^\tau}\\
	&\times \Arr{-\dfrac{\bfp(ci+t_{k+1})}{cd}+\bfp-1}{\bfp-1-\bfx_{\bfu_k,i}^\tau}\in \BZ.
	\end{split}
\end{equation}
Then
	\[\begin{split}
		&H_{\mu,c,\bfp,e,d}\\
	\equiv&\prod_{k=1}^b\prod_{n=0}^{d-2}\sum_{\tau\in S_n^\circ} \sgn(\tau)\prod_{i=1}^n \Arr{d-1}{d-1-\bfy_{\bfu_k,i}^\tau}
	\times (cd)^{\bfp-1-\bfx_{\bfu_k,i}^\tau}
	\\
	&\times \Arr{(ci+t_{k+1})\ell+\bfp-1}{\bfp-1-\bfx_{\bfu_k,i}^\tau}\\
	=&h_{u,e,d}\prod_{k=1}^b\prod_{n=0}^{d-2}\prod_{i=1}^n (d-1)! (cd)^{\bfp-1-\bfx_{\bfu_k,i}^\tau} \bigl((ci+t_{k+1})\ell+\bfp-1\bigr)!\mod p
	\end{split}\]
Note that $d-1,(ci+t_{k+1})\ell+\bfp-1<p$.
Thus 
	\[\NP_{u,m}(f)=\NP_{u,T}(f)=P_{u,e,d}\iff p\nmid H_{\mu,c,\bfp,e,d}\]
for $p>(d-e)(2d-1)$.
\end{proof}

\begin{proof}[Proof of Corollary~\ref{cor:cong}]
Since $p\nmid H_{\mu,c,\bfp,e,d}$, we have $H_{\mu,c,\bfp,e,d}\neq 0$.
Hence $p'\nmid H_{\mu,c,\bfp,e,d}$ for any $p'>H_{\mu,c,\bfp,e,d}$.
Note that
	\[\sum_{k=1}^b u_k=\frac{p-1}{c}\sum_{k=1}^b t_k,\]
thus $H_{[0,d],u}^\infty$ only depends on $\mu,c,\bfp,d$.
Since
	\[P_{u,e,d}(n)-H_{[0,d],u}^\infty(n)=\frac{d-e}{bd(p-1)}\sum_{k=1}^bC_{u_k,n-1}\le \frac{(d-e)\ov n(d-1)}{d(p-1)}\]
tends to zero as $p$ tends to infinity, the result then follows.
\end{proof}

\begin{example}
Assume that $p\equiv 1\bmod d$ and $d\mid u_k$ for all $k$. Write $p=dk+1$ and $t=u_k$. Then
	\[R_i:=R_{i,0}=\ov{e^{-1}i},\quad R_i:=r_{i,0}=\ov{-e^{-1}i},
	\quad \bfC_{t,n}=n,\quad S_n^\circ=\set{1}\]
and $x_{t,i}^1=\frac{(p-1)i+t}{d}, y_{t,i}^1=0$. Since
	\[h_{n,k}=\left(\prod_{i=0}^n\left(\frac{(p-1)i+u_k}{d}\right)!\right)^{-1}\in\BZ_p^\times,\]
we obtain that the Newton polygons coincide $H_{[0,d],u}^\infty$.
\end{example}

\section{The case \texorpdfstring{$e=d-1$}{e=d-1}}

If $pi-\tau(i)+t\notin d\BN+e\BN$ for some $i$, then $x_{t,i}^\tau<0$.
Set $1/k!=0$ for negative integer $k$.
Then
	\[h_{n,k}=\sum_{\tau\in S_{u_k,n}^\bullet}\sgn(\tau)\prod_{i=1}^n\frac{1}{x_{u_k,i}^\tau!y_{u_k,i}^\tau!},\]
where $S_{t,n}^\bullet$ the set of $\tau\in S_n^*$ such that the size of $\set{i\in I_n^*\mid R_{i,\alpha}+r_{\tau(i),\alpha}\ge d}$ is $C_{t,n,\alpha}$.

\begin{lemma}\label{eq:estimate_factorial}
Denote by $c(j)=\Arr{-\alpha j+\beta}{j}$.

(1) If $u_i=\alpha v_i+\beta$ for any $i$, then the matrix
\begin{equation}\label{eq:matrix transform}
	\bigl(\Arr{u_i}{j}\cdot \Arr{v_i+n}{n-j}\bigr)_{0\le j\le n}
	\implies \left(c(j) v_i^{n-j}\right)_{0\le j\le n}
\end{equation}
by third elementary column transformations.

(2) If $u_i\equiv \alpha v_i+\beta\bmod p$ for any $i$,
then \eqref{eq:matrix transform} holds by third elementary column transformations, modulo $p$.
\end{lemma}
\begin{proof}
(1) Write
	\[\Arr{\alpha x+\beta}{j}=\sum_{t=0}^{j}c_t(j)\cdot \Arr{x+j}{t},\]
then $c_0(j)=c(j)$ and
\begin{equation}\label{eq:expansion}
	\begin{split}
	 &\Arr{u_i}{j}\cdot \Arr{v_i+n}{n-j}\\
	=&\sum_{t=0}^j c_t(j)\cdot \Arr{v_i+j}{t}\cdot\Arr{v_i+n}{n-j}\\
	=&\sum_{t=0}^j c_t(j)\cdot \Arr{v_i+n}{n-j+t}.
	\end{split}
\end{equation}
Hence by third elementary column transformations,
	\[	\bigl(\Arr{u_i}{j}\cdot \Arr{v_i+n}{n-j}\bigr)
	\implies \bigl(c(j)\cdot \Arr{v_i+n}{n-j}\bigr)
	\implies \left(c(j)v_i^{n-j}\right).\]

(2) In this case,  \eqref{eq:expansion} holds modulo $p$. The result then follows easily.
\end{proof}

\begin{proof}[Proof of Theorem~\ref{thm:examples}]
Since $p>c(d^2-d+1)$, we have $p>(d-e)(2d-1)$.
Denote by $t=u_k$ and $t_k$ the minimal non-negative residue of $p^{-k} \mu$ modulo $c$.
Then $t=\frac{t_{k+1} p-t_k}{c}$.
If $c>1$, then $t\ge\frac{p-(c-1)}{c}\ge d(d-1)$ and $t<\frac{(c-1)p}{c}\le p-d(d-1)$. If $c=1$, then $t=0$.

Assume that $0\le n\le d-2$.
Denote by
	\[R_i=R_{i,t}=\ov{e^{-1}(pi+t)}=\ov{-pi-t}=-pi-t+\ell_i d\]
and
	\[r_i=r_{i,t}=\ov{-e^{-1}i}=\ov{i}.\]
Then
	\[\set{d-r_i\mid i\in I_n^*}=\set{d,d-1,\dots,d-n}.\]
We have
	\[\bfC_{t,n}=\#\set{i\in I_n^*\mid R_i\ge d-n}\]
and
	\[S_n^\bullet=\set{\tau\in S_n^*\mid R_i+\tau(i)\ge d\ \text{for}\ R_i\ge d-n}.\]
For $R_i<d-n$, we have $R_i+\tau(i)<d$ and
	\[x_{t,i}^\tau=pi+t-\ell_i e-\tau(i),\quad y_{t,i}^\tau=-pi-t+\ell_i d+\tau(i);\]
for $R_i\ge d-n$, we have $R_i+\tau(i)\ge d$ and
	\[x_{t,i}^\tau=pi+t-\ell_i e+e-\tau(i),\quad y_{t,i}^\tau=-pi-t+\ell_i d-d+\tau(i).\]

If $\tau\notin S_n^\bullet$, there is $i$ such that $y_{t,i}^\tau<0$ or $x_{t,i}^\tau<0$. 
Denote by
	\[(u_i,v_i)=\begin{cases}
		(pi+t-\ell_i e,-pi-t+\ell_i d),&\text{ if }R_i<d-n;\\
		(pi+t-\ell_i e+e,-pi-t+\ell_i d-d),&\text{ if }R_i\ge d-n.
	\end{cases}\]
Then
	\[h_{n,k}=\det\left(\frac{1}{(u_i-j)!(v_i+j)!}\right).\]
Apply Lemma~\ref{eq:estimate_factorial}(2) with $\alpha=-d^{-1}e,\beta=t(1-d^{-1}e)$, we obtain that
\[\begin{split}
&h_{n,k}\cdot\prod_{i=0}^n u_i!\cdot(v_i+n)!\\
\equiv &\prod_{j=0}^n \Arr{d^{-1}e(j-t)+t}{j} \cdot \det\left(v_i^{n-j}\right)\\
\equiv &\prod_{j=0}^n \Arr{d^{-1}e(j-t)+t}{j} \cdot \prod_{0\le i<j\le n}(v_i-v_j)\mod p.
\end{split}\]

If $R_i<d-n$, then $v_i=R_i\ge 0$; if $R_i\ge d-n$, then $v_i+n=R_i-d+n\ge 0$. Hence $0\le v_i+n\le d-1$ are different and $(v_i+n)!, (v_i-v_j)\in\BZ_p^\times$ if $i\neq j$.
Note that $u_i=\ell_i-R_i$ or $\ell_i-R_i+e$.
When $c=1$, we have $t=R_0=\ell_0$, $u_0=0$ or $e$, and for $i\ge 1$,
	\[u_i\ge \ell_i-R_i\ge \frac{pi+t}{d}-d+1\ge\frac{p}{d}-d+1\ge 0.\]
When $c>1$, we have
	\[u_i\ge \ell_i-R_i\ge \frac{pi+t}{d}-d+1\ge\frac{t}{d}-d+1\ge 0.\]
Meanwhile,
	\[u_i\le \ell_i-R_i+e=\frac{pi+t-(d-1)R_i+de}{d}\le\frac{p(d-2)+t+d e}{d}<p,\]
hence $u_i!\in\BZ_p^\times$.

For any $0\le k\le j-1$, we have
	\[0<e(j-t)+d(t-k)=d(j-k)+t-j\le (d-1)j+p-d(d-1)<p,\]
which means that $p\mid  \Arr{d^{-1}e(j-t)+t}{j}$.
Hence $h_{n,k}\in\BZ_p^\times$.
\end{proof}

\textbf{Acknowledgments.}
The author would like to thank Chuanze Niu and Daqing Wan for helpful discussions.
The author is partially supported by NSFC (Grant No. 12001510), Anhui Initiative in Quantum Information Technologies (Grant No. AHY150200) and the Fundamental Research Funds for the Central Universities (Grant No. WK0010000061).

\bibliographystyle{alpha}
\bibliography{mybib}
\end{document}